\documentclass[11pt,reqno]{amsart}
\usepackage{amssymb}
\usepackage{amsthm}
\usepackage[dvipsnames]{xcolor}
\usepackage{graphicx}
\usepackage[dvips]{psfrag}
\usepackage[active]{srcltx}
\usepackage{amsaddr}
\usepackage[shortlabels]{enumitem}

\usepackage{multirow}
\usepackage{verbatim}
\usepackage{color,graphicx}
\usepackage{mathtools}
\usepackage{xfrac}

\DeclareMathOperator\arctanh{arctanh}

\newcommand{\Proof}{\noindent\textbf{Proof.}\quad}

\newtheorem{theorem}{Theorem}[section]
\newtheorem{corollary}{Corollary}

\newtheorem{lemma}[theorem]{Lemma}
\newtheorem{proposition}{Proposition}

\theoremstyle{definition}

\newcommand{\R}{\mathbb{R}}
\newcommand{\C}{\mathbb{C}}

\newcommand{\N}{\mathbb{N}}

\newcommand{\X}{\mathcal{X}}

\advance\textwidth by +1.0in \advance\textheight by +0.5in
\advance\oddsidemargin by -0.5in \advance\evensidemargin by -0.5in
\advance\topmargin by -0.5in

\title[Integrability of a Family of Lotka--Volterra Three Species Biological System]
{Integrability of a Family of Lotka--Volterra Three Species Biological System}

\author[Aween Karim, Azad Amen, Waleed Aziz]
{Aween Karim$^{1}$,  Azad Amen$^{1,2,3,4}$ and Waleed Aziz$^5$}

\address{$^1$ Department of Mathematics, Collage of Basic Education, Salahaddin University--Erbil, Kurdistan region--Iraq}
\address{$^2$ Department of Mathematics, Basic Education College, Raparin University--Ranya, Kurdistan region--Iraq}
\address{$^3$ Department of Mathematics, Faculty of Science, Soran University--Soran, Kurdistan region--Iraq}
\address{$^4$ Department of Mathematics, Faculty of Science, Duhok University--Duhok, Kurdistan region--Iraq}
\address{$^5$ Department of Mathematics, College of Science, Salahaddin University--Erbil, Kurdistan region--Iraq}
\email{awen.karim@su.edu.krd}
\email{azad.amin@su.edu.krd}
\email{waleed.aziz@su.edu.krd}



\date{\today}

\begin{document}
\begin{abstract}
The aim of this study is to analyze the integrability problem of Lotka--Volterra three species biological
system. The system which considered in this work is a biological plausibility or a chemical model.
The system has a complex dynamical behavior because it is chaotic system. We, first show that the
system is a complete integrable when two of the involved parameters in the system are zero.
Second, thorough invariant algebraic surfaces and exponential factors, the nonintegrability problems have been investigated. Particularly, we show the non-existence of polynomial, rational, formal series, and Darboux first integrals when parameters are strictly positive.
\end{abstract}

\maketitle

\section{Introduction}
During the last fifty years there has been increasing interest in studying the autonomous differential system, mainly due to their many applications in natural science.
The conservation of ecological and biological systems is a primary concern for scientists and researchers, and controlling and analyzing the complex dynamical behavior of these systems is a challenge. The predation and competition species are the most popular interactions in which nonlinear polynomials are involved to represent such kind of interactions. The Lotka-Voltera system is one of the most prominent among existing these models. The authors in \cite{arneodo1980occurence} have been demonstrated that three species can possess chaotic behavior in an ecosystem. Samardzija and Greller \cite{samardzija1988explosive} derived a three species Lotka--Volterra biological system, which can be describe by a three dimensinal system of differential equations
\begin{equation}\label{system 1}
\begin{aligned}
\dot x &=x(1-y+cx-axz),\\
\dot y&=y(-1+x),\\
\dot z&=z(-b+ax^2).
\end{aligned}
\end{equation}   
In (\ref{system 1}), $x$ is the prey population and $y$, $z$ are predator populations, and by the biological meaning of the system, we only consider the parameters that are assumed to be non negative real numbers. For a certain value of parameters, $a=2.9851, b=3,$ and $c=2,$ the authors in \cite{samardzija1988explosive} showed that the three-species biological system (\ref{system 1}) was chaotic and it had very complicated dynamical behavior. In 1999, Costello \cite{costello1999synchronization} studied chaos synchronization in the integer-order system \eqref{system 1} for a certain parametric. Elsadany et al. \cite{elsadany2018dynamical} studied the dynamical behaviors of the system (1). They investigated the boundedness, existence, and uniqueness of the solutions of the system (1) and determined the stability and bifurcation of its equilibrium points.

One of the main difficulties in studying these differential systems consists of controlling the existence and nonexistence of first integrals and complete integrability.
Despite these earlier studies, the integrability of the system (\ref{system 1}) has not yet been investigated. The study of the integrability, the existence of first integrals,  is a main problem in the theory of differential equations. The existence of first integrals of a system are crucial in understanding the dynamics of the system. Therefore, it is important to understand whether a system has first integrals and whether they are analytical, smooth, etc. We study the existence of first integrals of system (\ref{system 1}), that described either by a formal power series (see \cite {bachmann1977henrici}, chapter one), rational function or by Darboux function using the Darboux theory of integrability (see \cite{darboux1878memoire}, \cite{christopher2000integrability} and \cite{zhang2017integrability}). Therefore, we are interested in studying the integrability of system \eqref{system 1} completely and  showing that it is non-integrable for most values of the parameters $a, b$ and $c$. \\
First, we look at the formal first integrals as it is a classical tool in studing differential equations. It is alos been used to described solutions around singularities \cite{seidenberg1968reduction}, the existence of first integrals given by formal series  \cite{nemytskii2015qualitative}, Moussu \cite{moussu1982symetrie}. The most success in using formal series to study differential equations has been achieved by Ecalle \cite{ecalle1992introduction}, who used them to prove the Dulac's conjecture. Mattei and Moussu \cite{mattei1980holonomie} are proved that formal integrability implies analytic integrability in dimension two.

The associated vector field of system \eqref{system 1}  is
\begin{equation}\label{eq vector}
\mathcal{X}=x(1-y+cx-axz)\frac{\partial}{\partial x}+y(-1+x)\frac{\partial}{\partial y}+z(-b+ax^2)\frac{\partial }{\partial z}.\\
\end{equation}

Let $U\subseteq \mathbb{C}^3$ be an open set. A non-constant function $H: U\to \mathbb{C}$ is a first integral of the polynomial vector field $\mathcal{X}$ on $U$  if it stays constant along all solution curves $(x(t), y(t), z(t))$ of (\ref{system 1}). Clearly $H$ is a first integral of $\mathcal{X}$ on $U$ if and only if $\mathcal{X}H=0$ on $U$. When $U=\mathbb{C}^3,$ the first integral $H$ is called a global first integral. When a first integral $H$ is a rational (polynomial or analytic) function, we say that $H$ is a rational (polynomial or analytic) first integral. when $H$ is a formal power series in the variables $x, y$ and $z$, we say that $H$ is a formal first integral.  Finally, a first integral is of Darboux type if it is of the form 
\begin{equation}\label{eq function}
f_{1}^{\lambda_{1}}\dots f_{p}^{\lambda_{p}} E_{1}^{\mu_{1}}\dots E_{q}^{\mu_{q}},
\end{equation}
where $f_{1},\dots,f_{p}$ are Darboux polynomials and $E_{1},\dots, E_{q}$ are exponential factors (see Section \ref{Preliminaries} for definitions), and $\lambda_{i}, \mu_{j}\in  \mathbb{C}$ for $i=1,\dots ,p,$ and $j=1, \dots,q.$
The functions of the form (\ref{eq function}) are called Darboux functions, and they are the base of the Darboux theory of integrability. 
The Darboux theory of integrability is an algebraic theory of integrability which based on the an adequate number of invariant algebraic surfaces and exponential factors associated to polynomial differential systems. In fact, to every Darboux polynomial there is associated some invariant algebraic surface, and the exponential factors appear when an invariant algebraic surface has multiplicity larger than $1$, see \cite{zhang2017integrability, christopher2007multiplicity, llibre2012darboux} for more details. Several well-known findings on the Darboux theory of integrability and analytical first integrals can be found in \cite{zhang2017integrability, aziz2021integrability, amenaintegrability}. The following is our first main results.

\begin{theorem}\label{Theorem invariant3D}
The unique irreducible Darboux polynomials of system \eqref{system 1} with non-zero cofactors are $x$ and $y$ and $z$.
\end{theorem}
This result is the basis of the Darboux theory of integrability, and its proof is given in Section \ref{sec:3D}.

\begin{theorem}\label{exponential}
The following statements hold for system \eqref{system 1}.
\begin{enumerate}[1)]
\item The unique exponential factors of system \eqref{system 1} are $\exp(x+z)$ and $\exp(y)$ with cofactors $cx^{2}-xy-bz+x$ and $y(x-1)$ respectively if $a>0$ and $c>0,$ with an exception if $c=0,$ an extra exponential factor $\exp((x+y+z)^{2})$ can be appear with cofactor $-2(x+y+z)(bz-x-y)$.

\item If $a=0, b>0$ and $c>0$, then system \eqref{system 1} has an exponential factor $\exp(z)$ with cofactor $-bz$  with an exception if $c=0,$ an extra exponential factor $\exp(x+y)$ can be derived with cofactor $x-y.$ Moreover, if $ab= 0,$ and $c>0,$, then it has no exponential factors.
\end{enumerate}
\end{theorem}

The proof of Theorem \ref{exponential} is given in Section \ref{sec:exponential factor}. The next two results state the existence of Darboux first integrals for a set of the values of the parameters and the absence of formal first integrals for system (\ref{system 1}).

\begin{theorem}\label{Darboux first integral}
If $a=c=0$, then for system \eqref{system 1} the following statements hold.\\
\begin{enumerate}[i)]
\item If $b=0,$ then it is Darboux integrable with the first integrals
\begin{center}$H_{1}=xy \exp(-x-y)  \quad \mbox{and} \quad    H_{2}=z.$\end{center}

\item If $b>0,$ then it is integrable with the following first integrals
\begin{center}
$H_{1}=xy \exp(-x-y)$  \quad and \quad  $H_{2}=z h(x,y),$\\
\end{center}
where $h(x,y)=\exp(\int ^{x}\!{\frac {b}{{\it s\left( {\rm LambertW} \left(-{\frac {xy{{\rm e}^{{\it s}-x-y}}}{{\it s}}}\right)+1\right)}}} d{\it s}),$ and lambertW computes the principal value of the Lambert W function.
\end{enumerate}
\end{theorem}
\begin{theorem}\label{has no Darboux type first integral}
System \eqref{system 1} with $a>0$ and $c>0$ has no first integrals of Darboux type.
\end{theorem}

\begin{theorem}\label{rational first integral}
If $a>0$ and $c>0$, then system \eqref{system 1} has no rational first integrals.
\end{theorem}

Darboux first integrals are not necessarily formal series first integrals, and vice versa. This demonstrates the independence of these topics. 
\begin{theorem}\label{formal first integral}
When $a> 0, c> 0$, system \eqref{system 1} does not admits any formal first integral.
\end{theorem}

To prove Theorem \ref{formal first integral}, when $b>0$, we rewrite system (\ref{system 1}) as a four dimensional system in the variables $x, y, z$, and $b$
\begin{equation}\label{system 7}
\dot x =x(1-y+cx-axz),\quad \dot y=y(-1+x),\quad  \dot z=z(-b+ax^2),\quad  \dot b=0.
\end{equation}

\begin{theorem}\label{formal in b}
Suppose that $a> 0,$ $c> 0,$ and $b> 0.$ If $f=f(x, y, z, b)$ is a formal first integral of system \eqref{system 7}, then $f$ is an arbitrary formal power series in the variable $b$.
\end{theorem}

From Theorem \ref{formal first integral}, directly the following result for the system \eqref{system 1} is obtained.
\begin{corollary}\label{corollary}
Suppose that $a> 0,$ and $c> 0.$ Then, the system \eqref{system 1} has no global analytic first integrals, no polynomial first integrals and no local analytic first integrals at the origin.
\end{corollary}

The proof of Theorems \ref {Darboux first integral}, \ref {formal first integral}, \ref{formal in b}, \ref{has no Darboux type first integral} and \ref{rational first integral} are given in Section \ref{sec:formal}.\\

Since the Darboux theory of integrability of a polynomial differential system is based on the existence of Darboux polynomials and their multiplicity, then the study of the existence or non-existence of Darboux first integrals needs to seek of the Darboux polynomials. We shall recall the main concepts of the Darboux theory of integrability. 

\section{Preliminaries}\label{Preliminaries}
A polynomial $f(x,y,z)\in \mathbb{C}[x,y,z]$ satisfying the equation
\begin{equation}\label{eq2}
x(1-y+cx-axz)\frac{\partial f}{\partial x}+y(-1+x)\frac{\partial f}{\partial y}+z(-b+ax^2)\frac{\partial f}{\partial z}=Kf,
\end{equation}
is called a Darboux polynomial of the system \eqref{system 1}, where $k=k(x,y,z)\in  \mathbb{C}[x, y,z],$ is a cofactor of $f(x,y,z)$ of degree at most two. Therefore, without loss of generality, the cofactor can be in the form
\begin{equation}\label{eq3}
K=\alpha_{0}+\alpha_{1}x+\alpha_{2}y+\alpha_{3}z+\alpha_{4}x^2+\alpha_{5}xy+\alpha_{6}xz+\alpha_{7}y^2+\alpha_{8}yz+\alpha_{9}z^2,
\end{equation}
where $\alpha_{i}\in \mathbb{C}$ for $i = 0,\dots,9$.
If $f (x,y,z)$ is a Darboux polynomial, then the surface $f(x,y,z)=0$ is an invariant algebraic surface.\\

We recall the following result in \cite{christopher2007multiplicity}.
\begin{lemma}\label{irreducible}
 Let $f$ be a polynomial and $f=\prod_{i=1}^s f_{j}^{\alpha_{j}}$ be its decomposition into irreducible factors in $\mathbb{C}[x, y, z].$ Then $f$ is a Darboux polynomial of system \eqref{system 1} if and only if all the $f_{j}$ are Darboux polynomials of system \eqref{system 1}. Moreover if $K$ and $K_{j}$ are the cofactors of $f$ and $f_{j},$ then $K=\sum{_{j=1}^{s}}\alpha_{j}K_{j}.$ 
\end{lemma}

An exponential factor $E$ of system (\ref{system 1}) is a function of the form $E=\exp(\frac{h}{f})$ where $h,f\in \mathbb{C}[x,y,z]$ satisfying $(h,f)=1$ and
\begin{equation}\label{eq3.1}
 \mathcal{X}E = LE,
\end{equation}
 for some polynomial $L= L(x, y,z)$ of degree at most two. Such function is called the cofactor of $E$.
 
\begin{proposition}\label{prop:exponential}
\cite{llibre2009darboux,llibre2009darboux1}. The following statements hold.
\begin{enumerate}[1)]
\item If $\exp({\frac{g}{h}})$ is an exponential factor for the polynomial differential system \eqref{system 1} and $h$ is not a constant polynomial, then $h=0$ is an invariant algebraic surface.
\item Eventually $\exp({g})$ can be an exponential factor, coming from the multiplicity of the infinite invariant plane.
\end{enumerate}
\end{proposition}

The following is well known result of the Darboux theory of integrability, for instance see Chapter 3 of \cite{zhang2017integrability}.
\begin{theorem}[Darboux Theory of Integrability]\label{Darboux Theory}
Suppose that a polynomial vector field $\mathcal{X}$ defined in $\mathbb{R}^{n}$ of degree $m$ admits $p$ Darboux polynomials $f_{i}$ with cofactors $K_{i}$ for $i=1,...,p,$ and $q$ exponential factors $E_{j}=\exp(g_{j}/h_{j})$ with cofactors $L_{j}$ for $j=1,...,q.$ If there exist $\lambda_{i}, \mu_{j}\in \mathbb{C}$ not all zero such that
\begin{equation}\label{eq condition}
\sum{_{i=1}^{p}}\lambda_{i}K_{i}+\sum{_{j=1}^{q}}\mu_{j}L_{j}=0,
\end{equation}
then the following real (multivalued) function is Darboux type
\begin{equation}\label{eq condition2}
f_{1}^{\lambda_{1}}...f_{p}^{\lambda_{p}} E_{1}^{\mu_{1}}...E_{q}^{\mu_{q}},
\end{equation}
and is also a first integral of the vector field $\X$.
\end{theorem}

The following result is well known, and it can be prove easily by Darboux theory of integrability.
\begin{lemma}\label{rational}
If system \eqref{system 1} has a rational first integral, then either it has a polynomial first integral or two Darboux polynomials with the same nonzero cofactor.
\end{lemma}

\section{Auxiliary Results}\label{sec:2D}
In this section, we state and prove some auxiliary results that will be used throughout the paper. These results study the existence of Darboux polynomials with a non-zero cofactor in two dimensions by restricting the original system to $y=0$ or $z=0$. It also examines the formal first integral of the system restricted to $z=0$.

\subsection{The system restricted to y=0}\label{subsec:y=0}
The system (\ref{system 1}) restricted to $y=0$ becomes
\begin{equation}\begin{split}\label{system 2}
&\dot x =x(1+cx-axz),\\
&\dot z=z(-b+ax^2).
\end{split}
\end{equation}
Let $\mathcal{Y}=\mathcal{X}|_{y=0}$ be the vector field associated system (\ref{system 2}),
\begin{equation*}
\mathcal{Y}=x(1+cx-axz)\frac{\partial }{\partial x}+z(-b+ax^2)\frac{\partial }{\partial z}.
\end{equation*}
Then the following results can be obtained.
\begin{lemma}\label{lemma(1)}
For system \eqref{system 2}, the following statements hold.
\begin{enumerate}[i)]
\item If $a>0,$ the unique irreducible Darboux polynomials with non-zero cofactors are $x$ and $z$.
\item If $a=0,$ the unique irreducible Darboux polynomials with non-zero cofactors are $x, z,$ and $1+cx$.
\end{enumerate}
\end{lemma}
\Proof
It can be easily verified that $x$ and $z$ are Darboux polynomials with respective cofactors, $1+cx-axz$ and $-b+ax^2$. First, suppose $a>0$, and we show that there is no other Darboux polynomial. Suppose that $f$ is an irreducible Darboux polynomial of system \eqref{system 2} with degree greater than or equal $2$ with a non-zero cofactor $k_{1}=\beta_{0}+\beta_{1}x+\beta_{2}z+\beta_{3}x^2+\beta_{4}xz+\beta_{5}z^2$. Then $f$ must stisfies
\begin{equation}\label{eq4}
x(1+cx-axz)\frac{\partial f}{\partial x}+z(-b+ax^2)\frac{\partial f}{\partial z}=k_{1}f.
\end{equation}
By restricting (\ref{eq4}) to the invariant plane $x = 0$ and denoting $f$ by $\bar{g}$, we obtain
\begin{equation}\label{eq5}
-bz\frac{d \bar{g}}{d z}=(\beta_{0}+\beta_{2}z+\beta_{5}z^2)\bar{g}.
\end{equation}
If $b=0$, we see that $(\beta_{0}+\beta_{2}z+\beta_{5}z^2)\bar{g}=0$ which implies $\beta_{0}=\beta_{2}=\beta_{5}=0$. When $b\ne0$, the solution of (\ref{eq5}) is
\[
\bar{g}={\it \ \tilde{d_{1}}}\,{z}^{-{\frac {{\it \beta_{0}}}{b}}}{\exp({\,{\frac {-z \left( {\it \beta_{5}}\,z+2\,{\it \beta_{2}} \right) }{2b}}}}),\quad \tilde{d_{1}}\in{\C}.
\]
Given that $\bar{g}$ is a polynomial, then must be $\beta_{2}=\beta_{5}=0$ and $\beta_{0}=-bm_{0},$ where $m_{0}\in \N\cup\left\{0\right\}$. Hence $k_{1}=-m_{0}b+\beta_{1}x+\beta_{3}x^2+\beta_{4}xz$ and $f=\tilde{d_{1}}z^{m_{0}}+xh(x,z)$, where $h(x, z)$ is a polynomial in the variables $x$ and $z$.
It is obvious that $\tilde{d_{1}}\ne0$ because $f$ is irreducible.\\
Now, restricting (\ref{eq4}) to the invariant plane $z=0$ and denoting $f$ by $\tilde{g}$, we obtain
\begin{equation}\label{eq6}
x(1+cx)\frac{d \tilde{g}}{d x}=(-m_{0}b+\beta_{1}x+\beta_{3}x^2)\tilde{g}.
\end{equation}
Here we consider two distinct cases.\\

\textbf{Case 1}: Suppose that $c>0$, then
$$
\tilde{g}={\it \tilde{d}_{2}}\, \left( 1+cx \right) ^{{\it m_{0}}\,b+{\frac {{\it \beta_{1}}}{c}}-{\frac {{\it \beta_{3}}}{{c}^{2}}}}{x}^{-{\it m_{0}}\,b}{\exp({{\frac {{\it \beta_{3}}\,x}{c}}}}),\quad \tilde{d_{2}}\in{\C}.
$$
Since $\tilde{g}$ is a polynomial, the must be $\beta_{3}=0$ and $\beta_{1}=c(m_{1}-m_{0}b)$, with $m_{1}\in \N\cup\left\{0\right\}$. Therefore,
\begin{equation}\label{eq7}
k_{1}=-m_{0}b+c(m_{1}-m_{0}b)x+\beta_{4}xz,
\end{equation}
and $f={\it \tilde{d}_{2}}\, \left( cx+1 \right) ^{m_{1}}{x}^{-{\it m_{0}}\,b}+zT(x,z)$, where $T(x, z)$ is a polynomial, and since $f$ is irreduced, then must $\tilde{d_{2}}\ne0$.\\
We write $f=\sum{_{j=0}^{n}}f_{j}{(x,z)}$, where each $f_{j}=f_{j}(x,z)$ denotes a homogeneous polynomial of degree $j$ in $x$ and $z$. Obviously, $f_{n}\ne0$. We can deduce from $(\ref{eq4})$ and $(\ref{eq7})$ that the terms of degree $n+2$ satisfy
\[
-ax^{2}z\frac{\partial f_{n}}{\partial x}+ax^{2}z\frac{\partial f_{n}}{\partial z}=\beta_{4}xzf_{n}.
\]
Solving this linear partial differential equation, we obtain $f_{n}=G_{n}(x+z)x^{-\frac{\beta_{4}}{a}}$. 
where $G_{n}(x+z)$ is an arbitrary function in terms of $x+z$. Since $f_{n}$ must be a homogeneous polynomial, we get that $\beta_{4}=-am_{2}$ where $m_{2}\in \N\cup\{0\}$ and $G_{n}(x+z)\in {\C[x, z]\backslash\{0\}}$. Note that $f_{n}\ne0$ implies that $G_{n}(x+z)\ne0$. Furthermore, because $f_{n}$ has degree $n$, we can write
\begin{equation}\label{eq8}
f_{n}=\tilde{d_{3}}(x+z)^{n-m_{2}}x^{m_{2}},\quad \tilde{d_{3}}\in\C\backslash \{0\},
\end{equation}
and the cofactor becomes $k_{1}=-m_{0}b+c(m_{1}-m_{0}b)x-am_{2}xz$.
Computing the terms of degree $n+1$ in $(\ref{eq4})$, we obtain
\[
cx^{2}\frac{\partial f_{n}}{\partial x}-ax^{2}z\frac{\partial f_{n-1}}{\partial x}+ax^{2}z\frac{\partial f_{n-1}}{\partial z}=c(m_{1}-m_{0}b)xf_{n}-am_{2}xzf_{n-1}
\]
Hence
\begin{eqnarray*}
f_{{n-1}}&=&{\frac{\tilde{d_{3}}c}{a}}{x}^{{m_{2}}(x+z)^{n-{m_{2}}-1}}((-{m_{0}}\,b+{m_{1}}-n) \ln(-z)+\ln(x)({m_{0}}\,b-{m_{1}}+{m_{2}}))\\
 &   &  + x^{m_{2}}{G_{n-1}}(x+z),
\end{eqnarray*}
where $G_{n-1}(x+z)$ is an arbitrary function in terms of $x+z$. Since $f_{n-1}$ must be a homogeneous polynomial, we obtain $G_{n-1}(x+z)\in \C[x, z]\backslash \{0\}$ and $-m_{0}b+m_{1}-n=0$, and $m_{0}b-m_{1}+m_{2}=0$. As a result, $n=m_{2}$, and therefore
\[
f_{n-1}=x^{n}G_{n-1}(x+z),
\]
which implies that $G_{n-1}(x+z)=0$ and hence $f_{n-1}=0$. Since  $\beta_{4}=-am_{2}=-a(m_{1}-m_{0}b)$, and computing the terms of degree $n-1$ in \eqref{eq4}, that satisfy
\begin{equation}\label{dn-1}
x\frac{\partial f_{n}}{\partial x}-ax^{2}z\frac{\partial f_{n-2}}{\partial x}+ax^{2}z\frac{\partial f_{n-2}}{\partial z}-bz\frac{\partial f_{n}}{\partial z}=-a(m_{1}-m_{0}b)xzf_{n-2}-bm_{0}f_{n}.
\end{equation}
The function
\[
f_{{n-2}} \left( x,z \right) ={\frac {{x}^{-{\it m_{0}}\,b+{\it m_{1}}}\tilde{d_{3}}{\it m_{1}}\, \left( \ln \left( x\right) x-x\ln  \left( -z \right) -x-z\right) }{a \left( x+z \right) ^{2}x}}+{x}^{-{\it m_{0}}\,b+{\it m_{1}}}{\it G_{n-2}} \left( x+z \right),
\]
satisfy \eqref{dn-1}
where $G_{n-2}(x+z)$ is an arbitrary function in terms of $x+z$. It obvious 
$\tilde{d_{3}}\ne0$, otherwise, $f$ becomes reducible. Since $f_{n-2}$ must be a homogeneous polynomial, we must have $m_{1}=0$.
Hence, $m_2=-m_{0}b$, and equation
\eqref{eq8} becomes
\[f_{n}=\tilde{d_{3}}x^{-m_{0}b},\quad \tilde{d_{3}}\in{\C\backslash \{0\}}.
\]
Since $f$ is a homogeneous polynomial and $m_{0}b\ge0$, $b\ne0$, then it must be $m_{0}=0$, and
hence $m_{2}=0$. As a result, $ \beta_{0}=\beta_{1}=\beta_{4}=0$, which implies that
the cofactor $k_1=0$.\\

\textbf{Case 2}: When $c=0$. In this case, equation (\ref{eq6}) becomes
\begin{equation}\label{eq9}
x\frac{\partial \tilde{g}}{\partial x}=(-m_{0}b+\beta_{1}x+\beta_{3}x^2)\tilde{g}.
\end{equation}
By solving it, we get $\tilde{g}={\it \tilde{d_{1}}}\,{x}^{-{\it m0}\,b}{\exp({\frac{1}{2}x \left( {\it \beta_{3}}\,x+2\,{\it \beta_{1}} \right) }})$, where $\tilde{d_{1}}\in{\C\backslash \{0\}}$.
We must have $\beta_{1}=\beta_{3}=0$, so
\begin{equation}\label{eq10}
k_{1}=-m_{0}b+\beta_{4}xz.
\end{equation}
When $c=0$, from (\ref{eq10}) and (\ref{eq4}), one can see
\begin{equation}\label{eq11}
x(1-axz)\frac{\partial f}{\partial x}+z(-b+ax^2)\frac{\partial f}{\partial z}=(-m_{0}b+\beta_{4}xz)f.
\end{equation}
Writing $f=\sum{_{j=0}^{n}}f_{j}{(x,z)}$ as a homogeneous polynomial and repeating the same way used in Case 1, we obtain
\begin{equation}\label{eq12}
f_{n}=\tilde{d_{2}}(x+z)^{n-m_{2}}x^{m_{2}},\qquad \tilde{d_{2}}\in{\C\backslash \{0\}},
\end{equation}
and
\[
f_{{n-1}}=x^{m_{2}}{G_{n-1}}\left(x+z\right),
\]
where $G_{n-1}(x+z)$ is a homogeneous polynomial in the terms of $(x+z)$. Then we can write
\[
f_{n-1}=\tilde{d_{3}}(x+z)^{n-m_{2}-1}x^{m_{2}},\qquad \tilde{d_{3}}\in{\C\backslash \{0\}}.
\]
Finally, we calculate the terms of degree $n$ in equation \eqref{eq11}, and we find that

\begin{eqnarray*}
f_{{n-2}}&=&{\frac{\tilde{d_{2}}}{a}}{x}^{-1+{m_{2}}}( x+z) ^{-2+n-{ m_{2}}}( -x ( {m_{0}}\,b+n) \ln( -z)+ x ( {m_{0}}\,b+n ) \ln ( x)\\
               &   & - ( x+z )( { m_{0}}\,b+{ m_{2}}))+ x^{ m_{2}}{G_{n-2}}( x+z ).
\end{eqnarray*}
Since $f_{n-2}$ must be a homogeneous polynomial and also since $\tilde{d_{2}}\ne0$, then must be
$m_{0}b+n=0$. We have $m_{0}b\ge0$, $b\ne0$ and $n$ is nonnegative, then must be $m_{0}=0$ and so $n=0$. If $n=0$, from (\ref{eq12}),
we get $f_{0}=\tilde{d_{2}}(x+z)^{-m_{2}}x^{m_{2}}$. Since the degree of $f_{0}$ is zero, it must
also $m_{2}=0$. Based on these values of parameters, we get that $\beta_{0}=\beta_{4}=0$,
and so $ k_{1} = 0$. This concludes the proof of statement (i).\\\\
Now considering that system (\ref{system 2}) with $a=0$ has the rational first integral $H=\frac{zx^b}{(1+cx)^b}$,
Proposition (\ref{rational}) implies that all invariant algebraic curves are in the set $\{(x, z)\in{\R^2} | Czx^b-(1+cx)^b=0,\: C\in \R\} \backslash \{x=-1/c\}$. This concludes the proof of the Lemma.\qed

\subsection{The system restricted to z=0}\label{subsec:z=0}
If system (\ref{system 1}) restricted to $z=0$, then it becomes
\begin{equation}\begin{split}
\dot x &=x(1-y+cx),\\
\dot y &=y(-1+x),\\
\end{split}
\label{system 3}
\end{equation}
The associated vector field of (\ref{system 3}) is
\[
\mathcal{Z}=\mathcal{X}|_{z=0}=x(1-y+cx)\frac{\partial }{\partial x}+y(-1+x)\frac{\partial }{\partial y}.
\]
To prove the next lemma, we need the following definition.
A polynomial $F(x, y, z)$ is said to be weight homogeneous of degree $r\in \N$ with respect to the weight exponent $s=(s_{1}, s_{2}, s_{3})$ for all $\mu \in \R\backslash {0}$ if it satisfies
\[
F(\mu^{s_{1}}x, \mu^{s_{2}}y, \mu^{s_{3}}z)=\mu^{r}F(x,y,z).
\]

\begin{lemma}\label{lemma(2)}
The unique irreducible Darboux polynomials of system \eqref{system 3} with non-zero cofactors are $x$ and $y$.
\end{lemma}
\Proof
Assume that $f$ is a Darboux polynomial of system (\ref{system 3}) with the cofactor $k_{2}$. Then
it must satisfies
\begin{equation}\label{eq13}
x(1-y+c x)\frac{\partial f}{\partial x}+y(-1+x)\frac{\partial f}{\partial y}=k_{2}f.
\end{equation}
It is clear that $k_{2}$ is a polynomial of degree at most one. Without loss of generality, we set $k_{2}=\alpha_{0}+\alpha_{1}x+\alpha_{2}y$. It is also obvious that $x$ and $y$ are Darboux polynomial
with respective cofactors $1-y+cx$ and $-1+x$. We now show that the system (\ref{system 3}) has no other Darboux polynomials. Suppose that $f$ is an irreducible Darboux polynomial of system (\ref{system 3}) of degree at least two, and we claim that $\alpha_{0}=\alpha_{1}=\alpha_{2}=0$.\\
On the invariant plane $x = 0$, the equation (\ref{eq13}) reduces to
\[
-y\frac{\partial \bar{h}}{\partial y}=(\alpha_{0}+\alpha_{2}y)\bar{h},
\]
where $\bar{h}=f|_{x=0}$, and its solution is
$\bar{h}=d_{1}{exp({-{\it \alpha_{2}}\,y}}){y}^{-{\it \alpha_{0}}},$ where $d_{1}\in{\C}$.
Since $\bar{h}$ must be a polynomial, then $\alpha_{2}$ must be zero, and $\alpha_{0}=-m_{0}$ where
$m_{0}\in \N\cup\left\{0\right\}$. It follows that $k_{2}=-m_{0}+\alpha_{1}x$ and
$f=d_{1}y^{-m_{0}}+x g_{1}(x, y)$, where $g_{1}(x, y)$ is a polynomial in the variables $x$ and $y$.
Note that $d_{1}\ne0$, otherwise $f$ would be reducible, which is a contradiction.\\
Now, if we set $\tilde{h}=f|_{y=0}$, the restriction (\ref{eq13}) to $y=0$ is
\[
x(1+c x)\frac{\partial \tilde{h}}{\partial x}=(-m_{0}+\alpha_{1}x)\tilde{h}.
\]
Let $c\ne0$, then $\tilde{h}=d_{2}\left( c x+1 \right) ^{{\frac {{\it \alpha_{1}}}{c}}+{\it m_{0}}}{x}^{-{\it m_{0}}}$, where $d_{2}\in{\C}$. 
Since $\tilde{h}$ must be a polynomial, we obtain
$\alpha_{1}=c(m_{1}-m_{0})$, where $m_{1}\in \N\cup\left\{0\right\}$.
So $k_{2}=-m_{0}+c (m_{1}-m_{0})x$ and $f=d_{2}x^{-m_{0}}(c x+1)^{{\frac {{\it \alpha_{1}}}{c}+{\it m_{0}}}}+y g_{2}(x, y)$, where $g_{2}(x, y)$ is a polynomial in the variables $x$ and $y$. Since $f$ is irreducible, then $d_{2}\ne0$.\\
We consider now the change of variables $x = X, y=\mu{Y}$ with $\mu\in{\C\backslash \{0\}}$.
Then system (\ref{system 3}) becomes
\begin{equation}\begin{split}
&\dot X =X(1-\mu{Y}+cX),\\
&\dot Y=Y(-1+X),
\end{split}
\label{system 4}
\end{equation}
where the prime denotes derivative with respect to the variable $T$. We apply the transformation
above and put $F(X, Y)=\mu^{n}f ( X, \mu{Y})$, where $n$ is the highest weight degree in the weight
homogeneous components of $f$ in $x$ and $y$ with the weight $(0,1)$ and $k_{2}=-m_{0}+c(m_{1}-m_{0})X$. Then $F$ satisfies
\[
\frac{dF}{dT}=\frac{d\mu^{n}f}{dt}=\mu^{n}\frac{df}{dt}=\mu^{n}k_{2}f=k_{2}F.
\]
Suppose $F= F_{0}+\mu F_{1}+ \mu^{2} F_{2}+\dots+\mu^{l}F_{l}$, where $F_{i}$ is a weight homogeneous polynomial in $X$ and $Y$ with the weight degree $n-i$ for $i= 0, 1,\dots,l$, and $n\ge l$. Clearly $f = F|_{\mu=1}$.
By definition of a Darboux polynomial, we have
\[
X(1-\mu{Y}+c X)\sum{_{i=0}^{n}}\mu^{i}\frac{\partial F_{i}}{\partial x}+Y(-1+X)\sum{_{i=0}^{n}}\mu^{i}\frac{\partial F_{i}}{\partial y}=(-m_{0}+c(m_{1}-m_{0})X)\sum{_{i=0}^{n}}\mu^{i}F_{i}.
\]
We claim that $m_{1}=m_{0}=0$. By Equating the terms with $\mu^{i}$, for $i=0,1$, we obtain
\begin{equation}\label{eq14}
\begin{split}
&L(F_{0})=(-m_{0}+c(m_{1}-m_{0})X)F_{0},\\
&L(F_{1})=(-m_{0}+c(m_{1}-m_{0})X)F_{1}+XY\frac{\partial F_{0}}{\partial x},
\end{split}
\end{equation}
where $L$ is a vector field
\[
L=X(1+c X)\frac{\partial }{\partial X}+Y(-1+X)\frac{\partial }{\partial Y}.
\]
Solving the first equation in (\ref{eq14}) and pulling back its solution to the second equation, we obtain
\[
F_{0} \left( X,Y \right)= {\it G_{0}} \left( YX \left( Xc+1 \right) ^{{\frac {-1-c}{c}}} \right)  \left( Xc+1 \right) ^{{\it m_{1}}}{X}^{-{\it m_{0}}},
\]
\begin{equation*}\begin{split}
F_{1}(X,Y)= \frac{(X c+1)^{\it m_{1}} X^{-\it m_{0}}}{\Gamma({\frac {2 c-1}{c}})\Gamma({\frac {c-1}{c}})}\Bigg(-D(G_{0})\left(Y X (X c+1)^{({\frac {c-1}{c}})}\right) X^3  \Gamma({\frac {c-1}{c}})\\
      \Gamma({\frac {2 c-1}{c}}) Y^2 (X c+1)^{({\frac {-2 c-2}{c}})} \mbox{hypergeom}([1,{\frac {2 c-1}{c}}],[2],-X c) c\\
     +M_{1} (X c+1)^{({\frac {-c-1}{c}})}+ \Gamma({\frac {c-1}{c}})\Big ({\frac{1}{2}} M_{2} D(G_{0}) \left(Y X (X c+1)^{({\frac {-c-1}{c}})}\right)\\
      Y^2 X (X c+1)^{({\frac {-2 c-2}{c}})}+G_{1}\left(Y X (X c+1)^{{\frac {-c-1}{c}}}\right) \Gamma({\frac {2 c-1}{c}})\Big)\Bigg),
\end{split}\end{equation*}
where 
\begin{equation*}\begin{split}
M_{1}=G_{0}\big(Y X \left(X c+1\right)^{\frac {-c-1}{c}}\big) Y \Gamma ({\frac {2 c-1}{c}}) \Bigg(c \Big(\Psi({\frac {2 c-1}{c}}) m_{0}+m_{0} \ln(X)+m_{0} \ln(c)\\
+ \mbox{hypergeom} \left( [1,1,{\frac {2 c-1}{c}}],[2,2],-X c\right) \left(m_{0}-m_{1}\right)X c
+ m_{0} (\Gamma-1)\Big) X \Gamma \left({\frac {2 c-1}{c}}\right)\\
-{\frac{1}{2}} m_{0}  \mbox{hypergeom} \left([1,1,{\frac {3 c-1}{c}}],[2,3],-X c\right) \Bigg) \Gamma \left({\frac {2 c-1}{c}}\right) X^{2} c^{2}\\
-\Gamma \left({\frac {c-1}{c}}\right) \Bigg(X c (m_{0}-m_{1}) \Psi({\frac {c-1}{c}})+X c  \left(m_{0}-m_{1}\right) \ln(X)\\
+X c (m_{0}-m_{1}) \ln(c)+c \Gamma  \left(m_{0}-m_{1}\right) X-m_{0}\Bigg),
\end{split}\end{equation*}
and
\begin{equation*}\begin{split}
M_{2}=\Bigg(\Big(2 X (c-1) \Psi ({\frac {2 c-1}{c}})+2 X (c-1) \ln(X)+2 X (c-1) \ln(c)-2\\
+2\Gamma (c-1) X \Big)\Gamma({\frac {2 c-1}{c}})+\Big( (-2 \Psi({\frac {3 c-1}{c}})-2 \ln(X)-2 \ln(c)\\
+ (-2  \mbox{hypergeom} ([1,1,{\frac {3 c-1}{c}}],[2,2],-X c ) c\\
+2 \mbox{hypergeom} ([1,1,{\frac {3 c-1}{c}}],[2,2],-X c)) X-2 \Gamma+2) \Gamma({\frac {3 c-1}{c}})\\
+\Gamma({\frac {4 c-1}{c}}) X c\: \mbox{hypergeom} ([1,1,{\frac {4 c-1}{c}}],[2,3],-X c)\Big) c X\Bigg).\\
\end{split}\end{equation*}

where $G_{0},G_{1}\in {\C[X,Y]}$ and have weight degrees $n$ and $n-1$, respectively. Since $F_{1}$ is a polynomial and $c\ne0$, we must have $m_{0}=0$ and $m_{1}-m_{0}=0$. This gives, $k_{2}=0$.\\
When $c=0$, the equation (\ref{eq13}) takes the form 
\begin{equation}\label{eq15}
x(1-y)\frac{\partial f}{\partial x}+y(-1+x)\frac{\partial f}{\partial y}=-m_{0}f.
\end{equation}
Writing $f=\sum{_{j=0}^{n}f_{j}}$, where each $f_{j}=f_{j}(x,y)$ represents  a homogeneous polynomial of degree $j$ in $x$ and $y$, in which $f_{n}\ne0$. By computing the terms of degree $n+1$ in (\ref{eq15}), we obtain
\[
-xy\frac{\partial f_{n}}{\partial x}+xy\frac{\partial f_{n}}{\partial y}=0,
\]
and whose solution is that is 
$$f_{n}(x,y)=T_{1}(x+y),$$
where $T_{1}(x+y)$ is a homogeneous polynomial of degree $n$. The terms of degree $n$ in (\ref{eq15}) satisfies
\[x\frac{\partial f_{n}}{\partial x}-xy\frac{\partial f_{n-1}}{\partial x}-y\frac{\partial f_{n}}{\partial y}+xy\frac{\partial f_{n-1}}{\partial y}=-m_{0}f_{n}.
\]
Solving the partial differential differential above, we get
\begin{equation}\begin{split}
f_{{n-1}} \left( x,y \right)= &-\frac {1}{x+y}\left( \ln\left( x \right)+\ln  \left( -y \right) \right)\left( x+y \right) \mbox {D} \left( T_{1}\right)  \left( x+y \right)\\
&-\ln  \left( -y \right) m_{{0}}T_{1} \left( x+y \right) + \left( x+y \right) T_{2} \left( x+y \right)\\
&+\ln\left( x \right) m_{{0}}T_{1}\left( x+y \right).\\
\end{split}\end{equation}
\\
Since $f_{{n-1}}$ must be a polynomial, we obtain $m_{0}=0$, which lead to $\alpha_{0}=0$, and so $k_{2}=0$. The proof now is complete.\qed

\begin{proposition}\label{prop:Formal}
The system \eqref{system 3} has no formal first integrals when $c>0$.
\end{proposition}
\Proof
We assume that $f(x,y)$ is a formal first integral of system (\ref{system 3}) with $c\ne0$. Then it satisfies
\begin{equation}\label{eq36}
x(1-y+cx)\frac{\partial f}{\partial x}+y(-1+x)\frac{\partial f}{\partial y}=0.
\end{equation}
We write $f$ as $f=\sum{_{k\ge 0}^{}}f_{k}(x)y^{k}$,  where each $f_{k}(x)$  is a formal  power series in the variable $x$. 
Denoting the restriction of $f$ to $y=0$ by $f_{0}= f_{0}(x)$ in equation \eqref{eq36}. Then
\[
x(1+cx)\frac{d f_{0}}{d x}=0,
\]
 and its solution is $f_{0}=d_{0}$, where $d_{0}$ is a constant. Then we can write 
 $f= f_{0}+yg(x,y)=d_{0}+yg(x,y)$ where $g=\sum{_{k\ge 0}^{}}f_{k+1}(x)y^{k}$. Now, the function $g$ must satisfies the equation
\begin{equation}\label{eq37}
x(1-y+cx)\frac{\partial g}{\partial x}+y(-1+x)\frac{\partial g}{\partial y}=-(-1+x)g.
\end{equation}
It is suffices to show that
\begin{equation}\label{eq38}
f_{k+1}(x)=0\quad   \mbox{for}\quad   k\ge 0.
\end{equation}
Since the restriction of $g$ to $y=0$ is $f_{1}=f_{1}(x)$, then from (\ref{eq37}), we obtain
\[
x(1+cx)\frac{d f_{1}}{d x}=-(-1+x)f_{1},
\]
and its solution is
\begin{equation}\label{eq39}
f_{1}=d_{1}x(1+cx)^{-1-\frac{1}{c}},
\end{equation}
where $d_{1}$ is a constant. By calculating the coefficient of  $y$ in (\ref{eq37}), we see
\[
x(1+cx)\frac{d f_{2}}{d x}-x\frac{d f_{1}}{d x}=-2(-1+x)f_{2},
\]
or equivalently,
\begin{equation}\label{eq40}
x(1+cx)\frac{d f_{2}}{d x}+d_{1}x(-1+x)(1+cx)^{-2-\frac{1}{c}}+2(-1+x)f_{2}=0.
\end{equation}
Multiplying (\ref{eq40}) by $(1+cx)^{2+\frac{1}{c}}$, gives
\begin{equation}\label{eq41}
x(1+cx)^{3+\frac{1}{c}}\frac{d f_{2}}{d x}+d_{1}x(-1+x)+2(-1+x)(1+cx)^{2+\frac{1}{c}}f_{2}=0.
\end{equation}
Evaluating (\ref{eq41}) on $x=-\frac{1}{c}$, we obtain that $-\frac{d_{1}(-1-\frac{1}{c})}{c}=0$, and since $c>0$, this implies $d_{1}=0$. Then, from (\ref{eq39}), we obtain $f_{1}=0$.
We now assume that (\ref{eq38}) is true for $k = 0,\dots,l-2$, and we will prove it is also true for
$k = l-1$. By the induction hypothesis, we have
\[
f_{1}=\dots=f_{l-2}=0\quad  \mbox{and}\quad  g=\sum{_{k\ge 0}^{}}f_{k+1}(x)y^{k}=y^{l-1}\sum{_{k\ge 0}^{}}f_{k+l}(x)y^{k}.
\]
By introducing $g$ in (\ref{eq37}) and determining the coefficient of $y^{l-1}$,  we see that $f_{l}$ satisfies
\[
x(1+cx)\frac{d f_{l}}{d x}+ (-1+x)(l-1)f_{l}=-2(-1+x)f_{l},
\]
which yields
\begin{equation}\label{eq42}
f_ {l}=d_{l}x^{l+1}(1+cx)^{-\frac{(1+c)(1+l)}{c}},
\end{equation}
where $d_{l}$ is a constant. By calculating the coefficient of $y^{l}$ in (\ref{eq37}), we obtain
\[
x(1+cx)\frac{d f_{l+1}}{d x}-x\frac{d f_{l}}{d x}+l(-1+x)f_{l+1}=-2(-1+x)f_{l+1},
\]
or equivalently,

\begin{equation}\begin{split}\label{eq43}
d_{{l}}{x}^{l+1} \left( -1+x \right)  \left( l+1 \right)\left( cx+1\right) ^{-{\frac { \left( 1+c \right)  \left( l+1 \right) }{c}}-1}+\left( c{x}^{2}+x \right) {\frac {\rm d}{{\rm d}x}}f_{{l+1}} \left( x\right)\\
+f_{{l+1}} \left( x \right)  \left( -1+x \right)  \left( l+2\right)=0.\\
\end{split}\end{equation}

Multiplying (\ref{eq43}) by $(1+cx)^{1+\frac{(1+c)(l+1)}{c}}$, implies

\begin{equation}\begin{split}\label{eq44}
\left(\left( c{x}^{2}+x \right) {\frac {\rm d}{{\rm d}x}}f_{{l+1}}\left(x \right) +f_{{l+1}} \left( x \right) \left( -1+x \right) \left( l+2 \right) \right) \left( cx+1 \right) ^{{\frac {lc+l+2\,c+1}{c}}}\\
+{x}^{l+1}d_{{l}} \left( -1+x \right) \left( l+1\right)=0.\\
\end{split}\end{equation}

Evaluating (\ref{eq44}) on $x=-\frac{1}{c}$, we obtain that $-\frac{d_{l}( l+1) (-\frac{1}{c})^{l}( -1-\frac{1}{c})}{c}=0$, and since $c>0$,
this implies $d_{l}=0$. Then, from (\ref{eq42}), we obtain $f_{l}=0$, for $k = l-1$. So we have by
induction that $f_{k}=0$ for all $k\ge 1$,  which yields $g=0$.
Consequently, $f=d_{0}$, and hence the system (\ref{system 3}) has no formal first integrals.\qed


\section{{\large Darboux polynomials with non-zero cofactors}}\label{sec:3D}
It is obvious that $x$,$y$, and $z$ are Darboux polynomials of \eqref{system 1}.
We show that system (\ref{system 1}) has no any Darboux polynomial of degrees greater
than one. The proof of Theorem \ref{Theorem invariant3D} will follows by the following Lemmas.\\\\
\textbf{Proof of Theorem \ref{Theorem invariant3D}}. We consider two cases.\\\\
\textbf{Case 1.} When $a=0$, the system (\ref{system 1}) reduces to
\begin{equation}\begin{split}
&\dot x =x(1-y+cx),\\
&\dot y=y(-1+x),\\
&\dot z=-bz,\\
\end{split}
\label{system 4}
\end{equation}
then any cofactor of \eqref{system 4} must be the form
\begin{equation}\label{eqq}
k=\alpha_{0}+\alpha_{1}x+\alpha_{2}y+\alpha_{3}z.
\end{equation}
\begin{lemma}\label{lemma(3)}
Let $f$ be an irreducible Darboux polynomial of degree greater than one with non-zero cofactor
\eqref{eqq}. Then $\alpha_{0}=\alpha_{1}=\alpha_{2}=\alpha_{3}=0$.
\end{lemma}
\Proof
Assume that $f$ is a Darboux polynomial of degree $n\ge2$ in system (\ref{system 4}) with a non-zero cofactor $k$. Then $f$ satisfies
\begin{equation}\label{eq20}
x(1-y+cx)\frac{\partial f}{\partial x}+y(-1+x)\frac{\partial f}{\partial y}-bz\frac{\partial f}{\partial z}=kf.
\end{equation}
Suppose that $\bar f$ is a restriction of $f$ to $x = 0$, then from equation (\ref{eq20}), $\bar f$ satisfies
\begin{equation}\label{eq21}
-y\frac{\partial \bar f}{\partial y}-bz\frac{\partial \bar f}{\partial z}=(\alpha_{0}+\alpha_{2}y+\alpha_{3}z)\bar f,\\
\end{equation}
and $\bar f$ satisfies
\begin{equation}\label{eq22}
\bar f={y}^{-{\it \alpha_{0}}}{\it \ F1} \left( z{y}^{-b}\right) {\exp({-{\it \alpha_{2}}\,y-{\frac {{\it \alpha_{3}}}{b}}z}}).\\
\end{equation}
Since $\bar f$ must be a polynomial, we get $\alpha_{2}=\alpha_{3}=0$. 
We now expand $f$ in powers of the variable $z$ so that $f=\sum{_{j=0}^{n}}f_{j}{z}^j$, where each
$f_{j}=f_{j}(x,y)$ is a polynomial in the variables $x$ and $y$. The Darbox polynomial of system (\ref{system 4}) with the restriction $z = 0$ is given by $f= f_{0}$, and according to the irreducibility condition for $f$, $f_{0}\neq 0$. Moreover, system (\ref{system 4}) restricted to $z=0$ yields system (\ref{system 3}), and so Lemma \ref{lemma(2)} guarantees that $\alpha_{0}=\alpha_{1}=0$. Therefore, $k=0$, which means that the system (\ref{system 4}) has no Darboux polynomials of greater than one with non-zero cofactor.\qed
\\

\textbf{Case 2.} $a>0$. To prove this case, we shall need the following Lemmas.
\begin{lemma}\label{lemma(4)}
Let $f$ be an irreducible Darboux polynomial of greater than one with non-zero cofactor as in \eqref{eq3}. Then $\alpha_{5}=\alpha_{7}=\alpha_{8}=\alpha_{9}=0$, and $\alpha_{4}=N_{4}a, \alpha_{6}=-N_{6}a$, where $N_{4}, N_{6}\in {\N\cup\{0\}}$.
\end{lemma}
\Proof
Let $f(x,y,z)=\sum{_{i=0}^{n}}f_{i}{(x,y,z)}$ be an irreducible Darboux polynomial of system (\ref{system 1}) where $n\ge2$, where each $f_{i}$ is a homogeneous polynomial of degree $i$ for $i= 0,1,\dots,n$. Clearly $f_{n}\ne 0$. Then $f$ satisfies
\begin{align*}
x(1-y+cx-axz)\sum{_{i=0}^{n}}\frac{\partial f_{i}}{\partial x}+y(-1+x)\sum{_{i=0}^{n}}\frac{\partial f_{i}}{\partial y}\\
+ z(-b+ax^2)\sum{_{i=0}^{n}}\frac{\partial f_{i}}{\partial z}=K\sum{_{i=0}^{n}}f_{i}.
\end{align*}
Computing the homogeneous terms of degree $n+2$, we get
\[
-ax^{2}z\frac{\partial f_{n}}{\partial x}+ax^{2}z\frac{\partial f_{n}}{\partial z}=(\alpha_{4}x^{2}+\alpha_{5}xy+\alpha_{6}xz+\alpha_{7}y^{2}+\alpha_{8}yz+\alpha_{9}z^{2})f_{n}.
\]
The solution this partial differential equation above is
\begin{align*}
f_{n}=T_{n}( y,x+z)( -z) ^{{\frac {{\it \alpha_{4}}}{a
}}+{\frac {{\it \alpha_{5}}\,y}{a \left( x+z \right) }}+{\frac {{\it \alpha_{7}}\,{y}^
{2}}{a \left( x+z \right) ^{2}}}}
 {x}^{-{\frac {{\it \alpha_{5}}\,y}{a \left( x+z \right) }}-{\frac {{\it \alpha_{6}}}{a}}-{\frac {{\it \alpha_{7}}\,{y}^{2}}{a\left( x+z \right) ^{2}}}+{\frac {{\it \alpha_{9}}}{a}}}\\
        {\exp({{\frac {{\it \alpha_{7}}\,{y}^{2}}{a \left( x+z \right) x}}+{\frac {{\it \alpha_{8}}\,y}{ax}}+{\frac { \left( x+z \right) {\it \alpha_{9}}}{ax}}}}),
\end{align*}
where $T_{n}$ is an arbitrary function of $y$ and $x+z$. Since $f_{n}$ is a homogeneous polynomial of degree $n$, we must have $\alpha_{5}=\alpha_{7}=\alpha_{8}=\alpha_{9}=0,$ and $\alpha_{4}=N_{4}a, \alpha_{6}=-N_{6}a$ where $N_{4},N_{6}\in {\N\cup\{0\}}$.\qed
\\
Therefore, 
\begin{equation}\label{eq23}
K=\alpha_{0}+\alpha_{1}x+\alpha_{2}y+\alpha_{3}z+N_{4}a x^{2}-N_{6}a xz.
\end{equation}

\begin{lemma}\label{lemma(5)}
Let $f$ be an irreducible Darboux polynomial of degree greater than one with non-zero cofactor $K$ as in \eqref{eq23}. Then $\alpha_{0}=\alpha_{1}=\alpha_{3}=N_{4}=N_{6}=0.$
\end{lemma}
\Proof 
Let $f\in {\C[x, y, z]}$ be an irreducible Darboux polynomial of system (\ref{system 1}) of degree $n\ge 2$. Writing $f$ in powers of the variable $y$, that is, $f=\sum{_{i=0}^{n}}f_{i}(x,z)y^{i}$, where $f_{i}=f_{i}(x,z)$ is a polynomial in variables $x$ and $z$ for each $i$. Given that $f$ is irreducible, we conclude that $f_{0}\ne 0$, otherwise, $f$ would be divisible by $y$, which is a contradiction. Moreover, $f_{0}$ is a Darboux polynomial of system (\ref{system 1}) restricted to $y = 0$. Then, $f_{0}$ satisfies
\begin{equation}\label{eq24}
x(1+cx-axz)\frac{\partial f_{0}}{\partial x}+z(-b+ax^{2})\frac{\partial f_{0}}{\partial z}=(\alpha_{0}+\alpha_{1}x+\alpha_{3}z+N_{4}a x^{2}-N_{6}a xz)f_{0}, 
\end{equation}
which is (\ref{eq4}) with $\beta_{0}=\alpha_{0},\beta_{1} =\alpha_{1}, \beta_{2} = \alpha_{3}, \beta_{3} =N_{4}a, \beta_{4} =-N_{6}a, \beta_{5}=0$. Proceeding exactly the same way as in Lemma \ref{lemma(1)}, we get that $\alpha_{0}=\alpha_{1}=\alpha_{3}=N_{4}=N_{6}=0$.\qed
\\\\
By Lemma \ref{lemma(5)}, we obtain that $K=\alpha_{2}y$.

\begin{lemma}\label{lemma(6)}
Let $f$ be an irreducible Darboux polynomial of greater than one with non-zero cofactor $K=\alpha_{2}y$. Then $\alpha_{2}=0$.
\end{lemma}
\Proof 
Let $f$ be an irreducible Darboux polynomial of system (\ref{system 1}) of degree $n\ge 2$.
Without loss of generality, let $f=\sum{_{i=0}^{n}}f_{i}{(y,z)x^{i}}$, where each $f_{i}=f_{i}(y,z)$ are polynomials in the variables $y,z$. It is clear that $f_{0}\ne 0$, otherwise, $f$ would be divisible by $x$, which is a contradiction. Since $f_{0}$ is independent of $x$, it is a Darboux polynomial of system (\ref{system 1}) restricted to $x=0$ that satisfies
\begin{equation}\label{eq26}
-y\frac{\partial f_{0}}{\partial y}-bz\frac{\partial f_{0}}{\partial z}=(\alpha_{2}y)f_{0}.
\end{equation}
We write $f_{0}$ as the sum of its homogeneous parts as $f_{0}=\sum{_{i=0}^{m}}f_{0,i}(y,z)$, where each $f_{0,i}=f_{0,i}(y,z)$ is a homogeneous polynomial in its variables of degree $i$ and $0\le m\le n$. Since $f_{0}\ne 0$, then $f_{0,m}\ne 0$. Computing the terms of degree $m+1$ in (\ref{eq26}), we have $(\alpha_{2}y)f_{0,m}=0$, and since $f_{0,m}\ne 0$, so must be $\alpha_{2}=0$. This completes the proof of the Theorem \ref{Theorem invariant3D}.\qed
\\

\section{Exponential Factors}\label{sec:exponential factor}

{\bf{Proof of Theorem \ref{exponential}}}. Let $E= \exp(\frac{h}{f})$ be an exponential factor of the system (\ref{system 1}) with cofactor $ L$ where $h, f\in{\C[x, y, z]}$ with $(h,f)=1$.
Then, from the definition of exponential factor and in view of Proposition \ref{prop:exponential}, we have either $f$ is a constant, in this case we can take $f=1$, or $f$ is a Darboux polynomial of system (\ref{system 1}), in this case from Theorem \ref{Theorem invariant3D}, $E$ can be of the form $E=\exp(\frac{h}{x^{s_{1}}y^{s_{2}}z^{s_{3}}})$ with $s_{1}, s_{2}$, and $s_{3}$ are non-negative integers and such that $h\in{\C[x, y, z]}$ is coprime with $x, y, z$. Obviously, when $s_{1}= s_{2}=s_{3}=0$, we are in the previous case.\\ 
We first prove that the system (\ref{system 1}) has no exponential factors of the form $\exp(\frac{h}{x^{s_{1}}y^{s_{2}}z^{s_{3}}})$. For this purpose, we assume that $\exp(\frac{h}{x^{s_{1}}y^{s_{2}}z^{s_{3}}})$ is an exponential factor of the system and that at least one of the $s_{1}, s_{2}$, or $s_{3}$ is positive, and then we get a contradiction. Clearly, by (\ref{eq3.1}), $h$ satisfies,
\begin{equation}\label{eq27}
\dot{x}\frac{\partial h}{\partial x}+\dot{y}\frac{\partial h}{\partial y}+\dot{z}\frac{\partial h}{\partial z}-h(\frac{\dot{x}s_{1}}{\partial x}+\frac{\dot{y}s_{2}}{\partial y}+\frac{\dot{z}s_{3}}{\partial z})=Lx^{s_{1}}y^{s_{2}}z^{s_{3}},
\end{equation}
We distinguish the following cases.\\

\textbf{Case 1}: When $s_{1}>{0}$. Evaluating (\ref{eq27}) on $x=0$ and setting $\bar{h}=h|_{x=0}$, then
$\bar{h}$ satisfies
\[
-y\frac{\partial \bar{h}}{\partial y}-bz\frac{\partial \bar{h}}{\partial z}=\bar{h}\big({s_{1}}(1-y)-{s_{2}}-b{s_{3}}\big).
\]
and clearly $\bar{h}\ne 0$ because $(h,x)=1$. The right hand side of the  equation above has one degree more than the left hand side, and this is a contradiction.\\

\textbf{Case 2}: When $s_{1}=0$ and $s_{2}>{0}$. On $y=0$ and setting $\tilde{h}=h|_{y=0}$, then
from  (\ref{eq27}) it must satisfies
\begin{equation}\label{eq28}
x(1+cx-axz)\frac{\partial \tilde{h}}{\partial x}+z(-b+ax^2)\frac{\partial \tilde{h}}{\partial z}=\tilde{h}\big({s_{2}}(-1+x)+{s_{3}}(-b+ax^2)\big).
\end{equation}
Clearly, $\tilde{h}\ne 0$, since $(h,y)=1$. Moreover, we note that $\tilde{h}$ is a Darboux polynomial of the vector field $\mathcal{Y}$. Then, from Lemma \ref{lemma(1)}, we can consider the following subcases below.\\

\textbf{Subcase 2.1}. If $a>0$. Then system (\ref{system 2}) has Darboux polynomials $x$ and $z$, so we can write $\tilde{h}=Gx^{m_{1}}z^{m_{2}}$ with $G$ is a constant and $m_{1}$, $m_{2}$ are non-negative integers such that $m_{1}+m_{2}>0$. 
We substitute $\tilde{h}$ in equation (\ref{eq28}), and we deduce
\[
-am_{1}xz+am_{2}x^{2}+cm_{1}x-bm_{1}x-bm_{2}+m_{1}=as_{3}x^{2}-bs_{3}+s_{2}x-s_{2}.
\]
Then clearly $am_{1}=0$, and then $m_{1}=0$. This implies that $\tilde{h}(x,z)=\tilde{h}(z)=Gz^{m_{2}}$. By substituting $\tilde{h}(z)$ in (\ref{eq28}), and simplifying it, we can derive
\begin{equation}\label{eq29}
(ax^{2}-b)m_{2}=as_{3}x^{2}-bs_{3}+s_{2}x-s_{2}.
\end{equation}

It is also easy to clear that $s_{2}=0$, which is a contradiction.\\

\textbf{Subcase 2.2}. If $a=0$. In this case, $\tilde{h}$ must be of the form $\tilde{h}=G(1+cx)^{bv}-x^{bv}z^{v}$, where $G$ is a constant and $v$ is a non-negative integer. The cofactor of $\tilde{h}$ is $vbc\,x$, and from (\ref{eq28}), we obtain
\[
{s_{2}}(-1+x)+{s_{3}}(-b)=vbc\,x,
\]
which implies $s_{2}=-bs_{3}\le 0$, and this a contradiction.\\

\textbf{Case 3}: When $s_{1}=s_{2}=0$ and $s_{3}>{0}$. From equation (\ref{eq27}) with $z=0$,
we get 
\begin{equation}\label{eq30}
x(1-y+cx)\frac{\partial \hat{h}}{\partial x}+y(-1+x)\frac{\partial \hat{h}}{\partial y}=\hat{h}{s_{3}}(-b+ax^2),
\end{equation}
where $\hat{h}=h|_{z=0}$ and $\hat{h}\ne{0}$ because $h$ is coprime with $z$. We emphasize that
$\hat{h}$ is a Darboux polynomial of the vector field $\mathcal{Z}$ with a non-zero cofactor.
Therefore, by Lemma \ref{lemma(2)}, we get that $\hat{h}=Gx^{m_{1}}y^{m_{2}}$ where $G$ as a
constant and $m_{1}$ and $m_{2}$ as non-negative integers such that $m_{1}+m_{2}>0$.
Then from (\ref{eq30}), we see that
\begin{equation}\label{eq31}
(c\,m_{1}+m_{2})x-m_{1}y+m_{1}-m_{2}=s_{3}(-b+ax^2).
\end{equation}
We consider two subcases depending on $a$ and $b$. Obviously, at least one of $a> 0$ or $b> 0$ is required for the existence of the Darboux polynomial $z$.\\

\textbf{Subcase 3.1}. If $a>0$. Computing the coefficient of $x^2$ in (\ref{eq31}), we get $as_{3}=0$. Hence, $s_{3}=0$, which is a contradiction.\\

\textbf{Subcase 3.2}. If $a={0}$ and $b>0$.  Then equation (\ref{eq31}) implies
\begin{equation}\label{eq31.1}
(c\,m_{1}+m_{2})x-m_{1}y+m_{1}-m_{2}=-b s_{3}.
\end{equation}
Computing the coefficients of $y$ and $x$  in (\ref{eq31.1}), we obtain $m_{1}=m_{2}=0$.
Consequently, $s_{3}=0$, which is a contradiction.\\

\textbf{Case 4}: $s_{1}=s_{2}=s_{3}={0}$. In this case, equation (\ref{eq27}) becomes
\begin{equation}\label{eq32}
x(1-y+cx-axz)\frac{\partial h}{\partial x}+y(-1+x)\frac{\partial h}{\partial y}+z(-b+ax^2)\frac{\partial h}{\partial z}=L,
\end{equation}
where $L=l_{0}+l_{1}x+l_{2}y+l_{3}z+l_{4}x^2+l_{5}xy+l_{6}xz+l_{7}y^2+l_{8}yz+l_{9}z^2$, with $l_{i}\in {\C}$ for $i =0,\dots,9$.
Here, we consider different subcases depending on $a$ and $c$.\\

\textbf{Subcase 4.1}. When $a>0$ and $c>0$. We write $h$ in the form $h=\sum{_{j=0}^{n}}h_{j}(x,y,z)$,
where each $h_{j}$ is a homogeneous polynomial of degree $j$ in the variables $x,y,z$.
Firstly, assume that $n\ge 3$. The terms of degree $n+2$ in (\ref{eq32}) are
\[
-ax^{2}z\frac{\partial h_{n}}{\partial x}+ax^{2}z\frac{\partial h_{n}}{\partial z}=0,
\]
whose solution is $h_{n}=W_{n}(y,x+z)$, where $W_{n}$ is an arbitrary $C^1$ function. Using the fact
that $h_{n}$ is a homogeneous polynomial of degree $n$, we can write
$h_{n}=\sum{_{j=0}^{n}}a_{j}y^{n-j}(x+z)^{j}$, where $a_{j}\in {\C}$.
The terms of degree $n+1$ in (\ref{eq32}) satisfy 
\[
x(-y+cx)\frac{\partial h_{n}}{\partial x}-ax^{2}z\frac{\partial h_{n-1}}{\partial x}+xy\frac{\partial h_{n}}{\partial y}+ax^{2}z\frac{\partial h_{n-1}}{\partial z}=0.
\]
or equivalently,
\begin{align*}
-ax^{2}z\frac{\partial h_{n-1}}{\partial x}+ax^{2}z\frac{\partial h_{n-1}}{\partial z}+(cx^{2}-xy)\sum{_{j=0}^{n}}ja_{j}y^{n-j}(x+z)^{j-1}\\
+xy\sum{_{j=0}^{n}}(n-j)a_{j}y^{n-j-1}(x+z)^{j}=0,
\end{align*}
The function
\begin{align*}
h_{n-1}=\frac{1}{2a(x+z)}(-4 \arctanh(\frac{x-z}{x+z})yB+A((2xc+2zc-4y) \arctanh(\frac{x-z}{x+z})\\
-\ln(-(x+z)xz)c(x+z)))+W_{n-1}(y, x+z),
\end{align*}
satisfy the partial differential equations above where 
$A=\sum{_{j=0}^{n}}ja_{j}y^{n-j}(x+z)^{j-1}$, $B=\sum{_{j=0}^{n}}a_{j}(-n+j)y^{n-j-1}(x+z)^{j}$, and $W_{n-1}$ is an arbitrary function of $y$ and $x+z$. Since $h_{n-1}$ must be a homogeneous polynomial of degree $n-1$, we must have $A=0$ and $B=0$, which implies that $a_{j}=0$ for $j=0,1,\dots, n$, and thus $h_{n}=\sum{_{j=0}^{n}}a_{j}y^{n-j}(x+z)^{j}=0$. This contradiction with the fact that $h_{n}$ is a polynomial of degree $n\ge{3}$. Then we must have $n\le{2}$. Therefore
\[
h=h_{0}+h_{1}x+h_{2}y+h_{3}z+h_{4}x^2+h_{5}xy+h_{6}xz+h_{7}y^2+h_{8}yz+h_{9}z^2.
\]
Substituting $h$ in (\ref{eq32}), we get that
\begin{align*}
l_{0}=0,  l_{1}=h_{1},  l_{2}=-h_{2},  l_{3}=-bh_{3},  l_{4}=ch_{1}+2h_{4}, l_{5}=-h_{1}+h_{2},  l_{6}= -bh_{6}+h_{6}, \\ l_{7}=-2h_{7},  l_{8}=-bh_{8}-h_{8}, l_{9}=-2bh_{9},
\end{align*}
where $h_{0}=h_{0}, h_{1}=h_{1}, h_{2}=h_{2}, h_{3}=h_{1}, h_{4}=h_{5}= h_{6}= h_{7}=h_{8}= h_{9}= 0$. Then $h=h_{0}+(x+z)h_{1}+yh_{2}$ with $L= (cx^2+(-y+1)x-bz)h_{1}+y(x-1)h_{2}$.\\\\

\textbf{Subcase 4.2}. When $a>0$ and $c=0$. From equation (\ref{eq32}), we have
\begin{equation}\label{eq33}
x(1-y-axz)\frac{\partial h}{\partial x}+y(-1+x)\frac{\partial h}{\partial y}+z(-b+ax^2)\frac{\partial h}{\partial z}=L.
\end{equation}
By proceeding in a similar way as above, we get
\[
h_{n-1}=\frac{y}{a(x+z)}(A+B)(\ln(x)-\ln(-z))+W_{n-1}(y, x+z),
\]
where $A=\sum{_{j=0}^{n}}ja_{j}y^{n-j}(x+z)^{j-1}, B=\sum{_{j=0}^{n}}a_{j}(-n+j)y^{n-j-1}(x+z)^{j}$, and $W_{n-1}$ is an arbitrary polynomial in $y$ and $x+z$. Since $h_{n-1}$ just admits polynomial solutions, we must have $A+B=0$, and $a_{j}=\frac{n!}{j! (n-j)!}A_{0}$ for $j=0,1,\dots,n$. Hence,
 $h_{n}=\sum{_{j=0}^{n}}\frac{n!}{j!(n-j)!}A_{0}y^{n-j}(x+z)^{j}$ and $h_{n-1}=W_{n-1}(y, x+z)$. Since $h_{n-1}$ must be a homogeneous polynomial of degree $n-1$, we write $h_{n-1}=\sum{_{j=0}^{n}}b_{j}y^{n-j-1}(x+z)^{j}$, where $b_{j}\in {\C}$.
Now, computing the homogeneous part of degree $n$ in (\ref{eq33}) yields
\[
-ax^{2}z\frac{\partial h_{n-2}}{\partial x}-xy\frac{\partial h_{n-1}}{\partial x}+x\frac{\partial h_{n}}{\partial x}-y\frac{\partial h_{n}}{\partial y}+xy\frac{\partial h_{n-1}}{\partial y}-bz\frac{\partial h_{n}}{\partial z}+ax^{2}z\frac{\partial h_{n-2}}{\partial z}=0,
\]
which is 

\begin{equation*}\begin{split}
-ax^{2}z\frac{\partial h_{n-2}}{\partial x}-xy\frac{\partial \sum{_{j=0}^{n}}b_{j}y^{n-j-1}(x+z)^{j}}{\partial x}+x\frac{\partial \sum{_{j=0}^{n}}\frac{n!}{j!(n-j)!}A_{0}y^{n-j}(x+z)^{j}}{\partial x}\\
-y\frac{\partial \sum{_{j=0}^{n}}\frac{n!}{j!(n-j)!}A_{0}y^{n-j}(x+z)^{j}}{\partial y}+xy\frac{\partial \sum{_{j=0}^{n}}b_{j}y^{n-j-1}(x+z)^{j}}{\partial y}\\
-bz\frac{\partial \sum{_{j=0}^{n}}\frac{n!}{j!(n-j)!}A_{0}y^{n-j}(x+z)^{j}}{\partial z}+ax^{2}z\frac{\partial h_{n-2}}{\partial z}=0.\\
\end{split}\end{equation*}

Solving it, we get

\begin{equation*}\begin{split}
h_{n-2}=\frac{1}{ax(x+z)^{2}}\big(-2 \arctanh(\frac{x-z}{x+z})xy(x+z)C-2 \arctanh(\frac{x-z}{x+z})xy(x+z)D\\
+ \left( -x \left( x-y+z \right) \ln  \left( -z \right) +x \left( x-y+z \right) \ln  \left( x \right)+\left( x+z \right) \left( bx+bz+y \right)  \right)\\
n A_{{0}} \left( x+z+y \right) ^{n-1}\big)+{\it W_{n-2}} \left( y,x+z \right),\\
\end{split}\end{equation*}

where $C=\sum{_{j=0}^{n}}(-n+j+1)b_{j}y^{n-j-2}(x+z)^{j}$, $D=\sum{_{j=0}^{n}}j b_{j}y^{n-j-1}(x+z)^{j-1}$,
and $W_{n-2}$ is an arbitrary polynomial in $y$ and $x+z$. Since $h_{n-2}$ must be a
homogeneous polynomial of degree $n-2$ and $n\ge 3$, then must be $A_{0}=0$, so $h_{n}=0$, 
which is a contradiction with the fact that $h_{n}$ is a polynomial of degree $n\ge{3}$. Then we must
have $n\le{2}$. Thus, $h=h_{0}+h_{1}x+h_{2}y+h_{3}z+h_{4}x^2+h_{5}xy+h_{6}xz+h_{7}y^2+h_{8}yz+h_{9}z^2$. By substituting $h$ in (\ref{eq33}), we can obtain
\[
 h=h_{0}+(x+y+z)^{2}h_{1}+(x+z)h_{2}+ yh_{3},
\]
and
\[
L=-2(x+y+z)(bz-x+y)h_{1}+\big((h_{3}-h_{2})y+h_{2}\big)x-bzh_{2}-yh_{3}.
\]
This completes the proof of the statement a.\\

\textbf{Subcase 4.3}. When $a={0}$ and $c>0$. In this case, we have
\begin{equation}\label{eq34}
x(1-y+cx)\frac{\partial h}{\partial x}+y(-1+x)\frac{\partial h}{\partial y}+z(-b)\frac{\partial h}{\partial z}=l_{0}+l_{1}x+l_{2}y+l_{3}z,
\end{equation}
where $l_{i}\in {\R}$ for $i =0,\dots,3$. First, we assume that $b>0$. We decompose $h$ as a sum of polynomials in the variable $x$ as $h=\sum{_{j=0}^{n}}h_{j}(y,z)x^{j}$, where each $h_{j}\in {\C[y,z]}$. Assume that $n\ge 1$, then the coefficient of $x^{n+1}$ in (\ref{eq34}), is
\[
c\, n\, h_{n}+y\frac{\partial h_{n}}{\partial y}=0,\]
and its solution is $h_{n}=W_{n}(z)y^{-c n}$.
Since $h_{n}$ must be a polynomial and $n\ge{1}$, $c>0$, then must be $ h_{n}=0$. As a result, $h_{n}=0$ if $n\ge{1}.$ For $n=0,$ we have $h=h_{0}(y,z),$ and the coefficients of $x$ from (\ref{eq34}) are
\begin{equation}\label{eq}
\begin{split}
y\frac{\partial h_{0}}{\partial y}&=l_{1},\\
-y\frac{\partial h_{0}}{\partial y}-bz\frac{\partial h_{0}}{\partial z}&=l_{0}+l_{2}y+l_{3}z.\\
\end{split}
\end{equation}  
The solution of the first equation in \eqref{eq}, is
\[
 h_{0}=l_{1}\ln(y)+W_{0}(z).
\]
The function $h_0$ is a polynomial if $l_{1}=0$, and putting it in the second equation in \eqref{eq},
we obtain
\[
W_{0}(z)=(-\frac{l_{2}}{b}y-\frac{l_{0}}{b})\ln(z)-\frac{l_{3}}{b}z.
\]
Since $W_{0}(z)$ must be a polynomial, we get $l_{0}=l_{2}=0$ and $W_{0}=-\frac{l_{3}}{b}z$. Then $h=-\frac{l_{3}}{b}z$ with $L=l_{3}z$.\\

If $b=0$, then equation (\ref{eq34}) becomes 
\begin{equation}\label{eq b0}
x(1-y+cx)\frac{\partial h}{\partial x}+y(-1+x)\frac{\partial h}{\partial y}=l_{0}+l_{1}x+l_{2}y,
\end{equation}
so by eliminating the variable $z$ and repeating the previous steps, we get
\[
c\, n\, h_{n}+y\frac{\partial h_{n}}{\partial y}=0,\quad \mbox{so that,}\quad h_{n}=d_{1}y^{-c n},\; d_{1}\in {\C}.
\]
Hence, $h_{n}=0$ if $n\ge{1}$. For $n=0$, we have $h=h_{0}(y)$, and by comparing the coefficients of $x$ from (\ref{eq b0}), we obtain
\begin{equation}
\begin{split}
y\frac{\partial h_{0}}{\partial y}&=l_{1},\\
-y\frac{\partial h_{0}}{\partial y}&=l_{0}+l_{2}y.\\
\end{split}
\end{equation}  
Solving the first equation, gives $h_{0}=l_{1}{\ln(y)}$, and it is a polynomial if $l_{1}=0$.
This implies that $h=h_{0}=0$. Eventually, the system (\ref{system 1}) has no exponential factors for $a=b=0$, and $c>0$.\\

\textbf{Subcase 4.4}. When $a={0}$ and $c=0$. The equation (\ref{eq32}) takes the form
\begin{equation}\label{eq35}
x(1-y)\frac{\partial h}{\partial x}+y(-1+x)\frac{\partial h}{\partial y}+z(-b)\frac{\partial h}{\partial z}=l_{0}+l_{1}x+l_{2}y+l_{3}z.
\end{equation}
If $b=0$, then the system has an exponential factor, as shown in the part (a) of Theorem \ref{Darboux first integral}. So from now on, we assume that $b\ne0$. We proceed in a similar way to the proof of Subcase 4.1. We write $h$ in the form $h=\sum{_{j=0}^{n}}h_{j}(x,y,z)$, where each $h_{j}$ is a homogeneous polynomial of degree $j$ in its variables. Assume $n\ge 2$. As before, we will show that $h$ has degree $n\le 1$. The terms of degree $n+1$ in (\ref{eq35}) are
\[
-xy\frac{\partial h_{n}}{\partial x}+xy\frac{\partial h_{n}}{\partial y}=0,
\]
and whose solution is, $h_{n}=W_{n}(x+y,z)$, where $W_{n}$ is an arbitrary polynomial in $x+y$ and $z$. Since $h_{n}$ is homogeneous of degree $n$, we can write $h_{n}=\sum{_{j=0}^{n}}a_{j}z^{n-j}(x+y)^{j}$, where $a_{j}\in {\C}$.  By computing the terms of degree $n$, in equation (\ref{eq35}), we can obtain
\[
h_{n-1}=-2\,b\,B\, arctanh(\frac{x-y}{x+y})-A\,\ln(-xy)+W_{n-1}(x+y,z),
\]
where $A=\sum{_{j=0}^{n}}ja_{j}(x+y)^{j-1}z^{n-j}$, $B=\sum{_{j=0}^{n}}a_{j}(n-j)(x+y)^{j-1}z^{n-j}$ and $W_{n-1}$ is an arbitrary function of $x+y$ and $z$. Since $h_{n-1}$ must be homogeneous polynomial of degree $n-1$ we must have $A=0$ and $B=0$ which implies that $a_{j}=0$ for $j=0,1,\dots,n$ and thus $h_{n}=0$, which is a contradiction. Then must be $n\le{1}$. Therefore $h=h_{0}+h_{1}x+h_{2}y+h_{3}z$.
Substituting $h$ in (\ref{eq35}) we get that $h=h_{2}(x+y)+h_{3}z$, and $L=h_{2}(x-y)-b\,h_{3}z$.
This completes the proof of the statement b.\qed


\section{First Integrals}\label{sec:formal}
The proof of Theorem $\ref{Darboux first integral}$ follows by direct computations, and it is omitted.\qed
\\

\textbf{Proof of Theorem \ref{has no Darboux type first integral}}.
According to Theorem \ref{Darboux Theory}, the system (\ref{system 1}) has a Darboux first integral if
and only if there exist $\lambda_{i},\mu_{j}\in {\C}$ not all zero, such that equation (\ref{eq condition})
is satisfied. It follows from Theorem \ref{Theorem invariant3D} that the system (\ref{system 1}) have
three Darboux polynomials with cofactors of $K_{1}=1-y+cx-axz, K_{2}=-1+x$, and $K_{3}=-b+ax^{2}$.
Now, by Theorem \ref{exponential}, when $a>0$ and $c>0$, there are two exponential factors
with cofactors 
$L_{1}=cx^{2}-xy-bz+x$ and $L_{2}=y(x-1)$. So equation (\ref{eq condition}) is equivalent to
\[
\lambda_{1}(1-y+cx-axz)+\lambda_{2}(-1+x)+\lambda_{3}(-b+ax^{2})+\mu_{1}(cx^{2}-xy-bz+x)+\mu_{2}(y(x-1))=0.
\]
Solving it, we get $\lambda_{1}=\lambda_{2}=\lambda_{3}= \mu_{1} =\mu_{2}= 0$. In short there are no first integrals of Darboux type in this case.\qed
\\

\textbf{Proof of Theorem \ref{rational first integral}}.
The proof of Theorem \ref{rational first integral} can be obtain easily from Theorem \ref{Theorem invariant3D} and corollary \ref{corollary} as well as Lemmas \ref{irreducible} and \ref{rational}.
\qed
\\

\textbf{Proof of Theorem \ref{formal first integral}}. If $b=0$, then system (\ref{system 1}) has the form
\begin{equation}\begin{split}
&\dot x =x(1-y+cx-axz),\\
&\dot y=y(-1+x),\\
&\dot z=ax^{2}z,\\
\end{split}
\label{system 6}
\end{equation}
Let $f$ be a formal first integral of system (\ref{system 6}). Without loss of generality, we can assume that $f$ has no constant terms. We write $f=\sum{_{j\ge 0}^{}}f_{j}(x,y)z^{j}$, where every $f_{j}(x,y)$  is a formal  power series in the variables $x$ and $y$. We consider two cases.

\textbf{Case 1}: If $f$ is not divisible by $z$. In this case, we have that $f_{0}=f_{0}(x,y)$ is
a formal first integral of system (\ref{system 6}) restricted to $z=0$. Note that system (\ref{system 6}), restricted to $z=0$, becomes system (\ref{system 3}). Consequently, $f_{0}$ is also a formal first integral of system (\ref{system 3}). However, we proved in Proposition \ref{prop:Formal} that system (\ref{system 3}) has no formal first integral. Therefore, $f_{0}$ is not a formal first integral of system (\ref{system 3}), and we have a contradiction.\\

\textbf{Case 2}: If $f$ is divisible by $z$. In this case, we write $f=z^{l}g$ where $l\ge1$, and $g$ is not divisible by $z$. Furthermore, $g$ is a formal power series that satisfy,
\[
x(1-y+cx-axz)\frac{\partial g}{\partial x}+y(-1+x)\frac{\partial g}{\partial y}+ax^{2}z\frac{\partial g}{\partial z}=a\,l\,x^{2}g.
\]
Now we introduce the variable $g=\exp(x)T$. Then, since $g$ is a formal series in the variables $x,y$, and $z$, we have that $T$ is a formal power series in the same variables, and it satisfies
\begin{equation}\label{eq45}
x(1-y+cx-axz)\frac{\partial T}{\partial x}+y(-1+x)\frac{\partial T}{\partial y}+ax^{2}z\frac{\partial T}{\partial z}=(-ax^{2}z+(al+c)x^{2}-xy+x)T,
\end{equation}
after dividing by $\exp(x)$. We write $T$ as $T=\sum{_{j\ge 0}^{}}T_{j}z^{j}$, where every $T_{j}=T_{j}(x,y)$ is a formal power series in the variables $x$ and $y$. Furthermore, since $g$ is not divisible by $z$, we have that $T$ is not divisible by $z$, so $T_{0}=T_{0}(x,y)\ne 0$. It follows that at least one of $\tilde{T}_{0}=\tilde{T}_{0}(x)\ne 0$ or $\bar{T}_{0}=\bar{T}_{0}(y)\ne 0$ must hold, where $\tilde{T}_{0}(x)$ is the restriction of $T_{0}$ to $y=0$, and $\bar{T}_{0}(y)$ is the restriction of $T_{0}$ to $x=0$. Otherwise, $T$ would be divisible by $z$, a contradiction. Without loss of generality, we can assume that $\tilde{T}_{0}\ne 0$. Moreover, if we restrict equation (\ref{eq45}) to $y=z=0$ and simplify it, we obtain that
\begin{equation}\label{eq46}
(1+cx)\frac{d \tilde{T}_{0}}{d x}=-(1+(al+c)x)\tilde{T}_{0}.
\end{equation}
Now we consider two cases.\\
\textbf{Subcase 2.1} $\tilde{T}_{0}$ is not divisible by $(1+cx)$. In this case, since $\tilde{T}_{0}\ne 0$ and $al\ne 0$ from (\ref{eq46}), we get a contradiction.\\
\textbf{Subcase 2.2} $\tilde{T}_{0}$ is divisible by $(1+cx)$. We write $\tilde{T}_{0}(x)=(1+cx)^{m}h(x)$ with $m\ge 1$ and $h(x)$ being a formal power series in the variable $x$ that is not divisible by $(1+cx)$. Then from (\ref{eq46}) we obtain $h$ must satisfy, after dividing by $(1+cx)^{m}$,
\begin{equation}\label{eq47}
(1+cx)\frac{d h}{d x}=-(1+(al+c)x-cm)h.
\end{equation}
Since $c m>0$  and $al>0$, we have from (\ref{eq47}) that $h$ must be divisible by $(1+cx)$, a contradiction. So system (\ref{system 1}) has no formal first integral if $b=0$. \qed\\

Now we consider system (\ref{system 1}) with $a\ne 0$,$c\ne 0,$ and $b\ne0$. We provide the proof  of Theorem \ref{formal in b} in order to prove Theorem \ref{formal first integral} when $b>0$.\\

\textbf{Proof of Theorem \ref{formal in b}}.
 Since $b$ is one of the parameters in the system, we can rewrite it as
\begin{equation}\begin{split}
&\dot x =x(1-y+cx-axz),\\
&\dot y=y(-1+x),\\
&\dot z=z(-b+ax^2),\\
&\dot b=0.\\
\end{split}
\label{system 8}
\end{equation}   
In other words, we add a new variable, $b$, that was a parameter in system (\ref{system 1}). Note that a non constant function $f = f(b)$ is a first integral of system  (\ref{system 8}), but it is not a first integral of the system (\ref{system 1}).\\
We assume that $f=f (x,y,z,b)$ is a formal power series first integral of system (\ref{system 8}). Expanding $f$ in powers of the variable $b$, we get $f=\sum{_{k\ge 0}^{}}f_{k}(x,y,z)b^{k}$, where each $f_{k}$ is a formal series in its variables. We can write $f=f_{0}+bg$ where $g=\sum{_{k\ge 0}^{}}f_{k+1}b^{k}=\sum{_{k\ge 0}^{}}f_{k+1}(x,y,z)b^{k}$ is a formal series in variables $x,y$, and $z$. Since $f(x,y,z,0)$ is a formal first integral of the system (\ref{system 1}) with $b=0$, and since $a\ne 0$, and $c\ne 0$, we are in the assumptions of Theorem \ref{formal first integral} for $b>0$, applying it, we get that $f(x,y,z,0)=f_{0}=d_{0}$. We claim that
\begin{equation}\label{eq48}
 f_{k+1}=d_{k+1},\quad  k\ge 0,
\end{equation}
where $d_{k+1}$ are constants for $k\ge 0$. Then, since $f$ is a first integral, it satisfies ${\X}f=0$. So, the function $g$ satisfies the equation
\begin{equation}\label{eq49}
x(1-y+cx-axz)\frac{\partial g}{\partial x}+y(-1+x)\frac{\partial g}{\partial y}+(b+ax^{2}z)\frac{\partial g}{\partial z}=0.
\end{equation}
Moreover, $f_{1}=g(x, y, z, 0)$ satisfies (\ref{eq49}) restricted to $b=0$, that is,
\[
x(1-y+cx-axz)\frac{\partial f_{1}}{\partial x}+y(-1+x)\frac{\partial f_{1}}{\partial y}+(ax^{2}z)\frac{\partial f_{1}}{\partial z}=0.
\]
Hence, $f_{1}$ is a formal power series first integral of system (\ref{system 1}) with $b=0$, and by assumption, additionally, we have that $a\ne 0$, and $c\ne 0$. So, from Theorem \ref{formal first integral} for $b>0$, we obtain that $f_{1}= d_{1}$, a constant. This proves (\ref{eq48}) for $k=0$. Now we assume (\ref{eq48}) is true for $k=0,\dots,j-1$ with $j\ge1$, and we will prove it for $k=j$. Since $g=bd_{1}+b\sum{_{k\ge 1}^{}}f_{k+1}b^{k}$, then $f=d_{0}+bd_{1}+b\sum{_{k\ge 1}^{}}f_{k+1}b^{k}$, and by the induction hypothesis, we have
\[
f=\sum{_{k=0}^{j-1}}d_{k}b^{k}+b^{j+1}\sum{_{k\ge j}^{}}f_{k+1}b^{k-j},
\]
so $\sum{_{k\ge j}^{}}f_{k+1}b^{k-j}$ is a formal first integral of system (\ref{system 7}). Consequently, $f_{j+1}$ is a formal first integral of system (\ref{system 1}) with $b = 0$, and from part (a), with $a\ne 0$,$c\ne 0$, we get $f_{j+1}=d_{j+1}$, and this proves the claim for $k=j$. Then, from (\ref{eq48}), we get that $f=\sum{_{k\ge 0}^{}}d_{k+1}b^{k}$, which finishes the proof of the theorem.\qed

\end{document}